\documentclass[11pt,reqno]{amsart}
\usepackage{graphicx}
\usepackage{verbatim}
\usepackage{textcomp}
\usepackage{amssymb}
\usepackage{cite}
\usepackage{amsmath}
\usepackage{latexsym}
\usepackage{amscd}
\usepackage{amsthm}
\usepackage{mathrsfs}
\usepackage{xypic}
\usepackage{bm}
\usepackage{url}
\usepackage{hyperref}

\newtheorem{thm}{Theorem}[section]

\newtheorem{lem}[thm]{Lemma}

\theoremstyle{definition}
\newtheorem{defn}{Definition}[section]

\theoremstyle{remark}
\newtheorem{rem}{Remark}[section]
\numberwithin{equation}{section}
\begin{document}
\title{Rigid properties of generalized $\tau$-quasi Ricci-harmonic metrics}
\author{Fanqi Zeng}
\address{School of Mathematics and Statistics, Xinyang Normal University, Xinyang, 464000, P.R. China} \email{fanzeng10@126.com}
\thanks{This work was supported by Nanhu Scholars Program for Young Scholars of XYNU}
\maketitle

\begin{abstract}
In this paper, we study compact generalized $\tau$-quasi Ricci-harmonic metrics. In the
first part, we explore conditions under which generalized $\tau$-quasi Ricci-harmonic metrics are harmonic-Einstein and give some characterization results for it. In the second part, we obtain some
rigidity results for compact $(\tau, \rho)$-quasi Ricci-harmonic metrics which are special case of generalized $\tau$-quasi Ricci-harmonic metrics. In the third part, we shall give
two gap theorems for compact $\tau$-quasi Ricci-harmonic metrics by showing
some necessary and sufficient conditions for the metrics to be harmonic-Einstein.
\end{abstract}

\section{Introduction}

In this paper, we investigate the rigid properties for generalized $\tau$-quasi Ricci-harmonic metrics and a special class of generalized $\tau$-quasi Ricci-harmonic metrics. To be precise, let us first give our notation. Throughout
this paper, let $(M, g)$ and $(N, h)$ be two static complete Riemannian manifolds of dimension $m$ and $n$, respectively. We denote by $Ric$ and
$R$ the Ricci tensor (with respect to $g$) and scalar curvature of $M$, respectively. We denote by $\nabla$, $\Delta$ and $\nabla^{2}$ the gradient, the Laplacian and the
Hessian on $(M, g)$, respectively. Let $\otimes$, $\langle, \rangle$ and $|\cdot|$ be the tensorial product, the metric $g$ and its associated norm, respectively.
Let $\phi: (M, g) \to (N, h)$ be a smooth map between $(M, g)$ and $(N, h)$, $f : M \to \mathbb{R}$ a smooth function on $M$.

First we give the precise definition of a generalized $\tau$-quasi Ricci-harmonic metric.
\begin{defn}\label{def1}
For $\tau>0$, we call a metric $g$ of $M$ \textbf{generalized $\tau$-quasi Ricci-harmonic} (with respect to $h$), if for some map $\phi: (M, g) \to (N, h)$, some potential function $f : M \to \mathbb{R}$ and a soliton function $\lambda: M \to \mathbb{R}$ and some constant $\alpha\geq 0$, $g$ satisfies the following coupled system
\begin{equation}\label{general1}
\left\{\begin{array}{l}
Ric_{f, \tau}-\alpha\nabla\phi\otimes\nabla\phi=\lambda g,\\
\tau_{g}\phi=\langle\nabla\phi,  \nabla f\rangle,
\end{array}\right.
\end{equation}
where $\nabla\phi\otimes\nabla\phi:= \phi^{\ast}h$ is the pull-back of the metric $h$ via $\phi$ and $\tau_{g}\phi=trace \nabla d\phi$ denotes the tension field of $\phi$ with respect to $g$ \cite{Eells}, and $Ric_{f, \tau}$ denotes the $\tau$-Bakry-\'{E}mery Ricci tensor
$$Ric_{f, \tau}=Ric+ \nabla^{2} f-\frac{1}{\tau}df\otimes df.$$
\end{defn}
For simplicity, we put
$$Ric_{\phi}:=Ric-\alpha \nabla \phi\otimes \nabla \phi$$
and call it a \textbf{generalized Ricci curvature}. We denote by
$$R_{\phi}:=trace(Ric_{\phi})=R-\alpha |\nabla \phi|^{2}$$
and call it a \textbf{generalized scalar curvature}.

Note that, if $(N, h)=(\mathbb{R}, dr^{2})$ and $\phi: (M, g) \to (\mathbb{R}, dr^{2})$ is a constant function in \eqref{general1},
then the generalized $\tau$-quasi Ricci-harmonic metric is exactly a generalized $\tau$-quasi Einstein metric, see \cite{Barros,Catino,Deng,Hu,Hu2017,Huang,Wanglf2012}.

In the special case $\lambda=\rho R_{\phi}+\mu$ in \eqref{general1} with $\rho,\mu$ two real constants, $g$ becomes a $(\tau, \rho)$-quasi Ricci-harmonic metric, which is defined as follows.
\begin{defn}\label{def2}
For $\tau>0$, we call a metric $g$ of $M$ \textbf{$(\tau, \rho)$-quasi Ricci-harmonic} (with respect to $h$), if there exists some map $\phi$, some potential function $f$ on $(M, g)$ and three real constants $\alpha\geq 0, \rho, \mu$ such that $g$ satisfies the following coupled system
\begin{equation}\label{general2}
\left\{\begin{array}{l}
Ric_{f, \tau}-\alpha\nabla\phi\otimes\nabla\phi=(\rho R_{\phi}+\mu) g,\\
\tau_{g}\phi=\langle\nabla\phi,  \nabla f\rangle.
\end{array}\right.
\end{equation}
\end{defn}
Note that, if $(N, h)=(\mathbb{R}, dr^{2})$ and $\phi: (M, g) \to (\mathbb{R}, dr^{2})$ is a constant function in \eqref{general2},
then the $(\tau, \rho)$-quasi Ricci-harmonic metric is exactly a $(\tau, \rho)$-quasi Einstein metric, see \cite{Huang13,Huang17,Shin,Wanglf2018}.

In the special case $\lambda$ is a constant in \eqref{general1}, $g$ becomes a $\tau$-quasi Ricci-harmonic metric, which was given in \cite{Wanglf,Wanglf2017}.
\begin{defn}\label{def3}
For $\tau>0$, we call a metric $g$ of $M$ \textbf{$\tau$-quasi Ricci-harmonic} (with respect to $h$), if for some map $\phi$, some potential function $f$ and a soliton constant $\lambda$ and some constant $\alpha\geq 0$, $g$ satisfies the following coupled system
\begin{equation}\label{general3}
\left\{\begin{array}{l}
Ric_{f, \tau}-\alpha\nabla\phi\otimes\nabla\phi=\lambda g,\\
\tau_{g}\phi=\langle\nabla\phi, \nabla f\rangle.
\end{array}\right.
\end{equation}
\end{defn}
Note that, if $(N, h)=(\mathbb{R}, dr^{2})$ and $\phi: (M, g) \to (\mathbb{R}, dr^{2})$ is a constant function in \eqref{general3},
then the $\tau$-quasi Ricci-harmonic metric is exactly a $\tau$-quasi Einstein metric, see \cite{Case,He,Wanglinf2018}. We know from \cite{Besse} that $\tau$-quasi Einstein metrics are closely
relative to the existence of warped product Einstein manifolds for any positive integer $\tau$.

It is important to point out that if $\tau=\infty$, the equation \eqref{general3} becomes gradient Ricci-harmonic soliton metric,
which was introduced by M\"{u}ller \cite{Muller}.
\begin{defn}\label{def4}
We call a metric $g$ of $M$ \textbf{gradient Ricci-harmonic soliton metric} (with respect to $h$), if for some map $\phi$, some potential function $f$ and a soliton constant $\lambda$ and some constant $\alpha\geq 0$, $g$ satisfies the following coupled system
\begin{equation}\label{general4}
\left\{\begin{array}{l}
Ric_{f}-\alpha\nabla\phi\otimes\nabla\phi=\lambda g,\\
\tau_{g}\phi=\langle\nabla\phi, \nabla f\rangle,
\end{array}\right.
\end{equation}
where $Ric_{f}$ denotes the $\infty$-Bakry-\'{E}mery curvature defined by:
$$Ric_{f}=Ric+ \nabla^{2} f.$$
We call a gradient Ricci-harmonic soliton metric shrinking, steady or expanding, if
$\lambda >0$, $\lambda =0$ or $\lambda <0$, respectively.
\end{defn}
As pointed out in \cite{Muller,Muller2012}, a gradient Ricci-harmonic soliton arises from the the Ricci-harmonic flow.
Let $(M, g(t))$ be a family of complete Riemannian manifolds with Riemannian metrics $g(t)$
evolving by the Ricci-harmonic flow
\begin{equation}\label{genera4}
\left\{\begin{array}{l}
\frac{\partial g}{\partial t}(t)=-2Ric_{g(t)}+2\alpha(t)\nabla\phi(t)\otimes\nabla\phi(t),~~g(0)=g_{0},\\
\frac{\partial \phi}{\partial t}(t)=\tau_{g(t)}\phi(t),~~\phi(0)=\phi_{0},
\end{array}\right.
\end{equation}
where $\alpha(t) \geq 0$ is a non-negative time-dependent coupling constant, $\phi(t): (M, g(t)) \to (N, h)$ is a family of smooth maps between $(M, g(t))$ and a fixed complete Riemannian manifold $(N, h)$ and $\nabla\phi(t)\otimes\nabla\phi(t):= \phi(t)^{\ast}h$ is the pull-back of the metric $h$ via $\phi(t)$. And as before, $\tau_{g}\phi=trace \nabla d\phi$ denotes the tension field of $\phi$ with respect to $g(t)$.

Note that, if $(N, h)=(\mathbb{R}, dr^{2})$ and $\phi: (M, g) \to (\mathbb{R}, dr^{2})$ is a constant function in \eqref{general4},
then the gradient Ricci-harmonic soliton metric is exactly a gradient Ricci soliton metric. The gradient Ricci soliton metrics play a very important
role in Hamilton's Ricci flow as they correspond to the self-similar solutions and often arise as singularity models, for a survey in
this subject we refer to the work due to Cao in \cite{Cao09}.

When $f$ is constant, the $\tau$-quasi Ricci-harmonic metric and gradient Ricci-harmonic soliton metric defined in \eqref{general3} and \eqref{general4} are called
\textbf{harmonic-Einstein}, which satisfy the following coupled system \cite{Williams}:
\begin{equation}\label{general5}
\left\{\begin{array}{l}
Ric-\alpha\nabla\phi\otimes\nabla\phi=\lambda g,\\
\tau_{g}\phi=0.
\end{array}\right.
\end{equation}
We point out that if $\lambda=\lambda(x)$ is a smooth function on $(M, g)$, then gradient Ricci-harmonic soliton and harmonic-Einstein are almost Ricci-harmonic soliton and almost harmonic-Einstein, respectively \cite{Abolarinwa19}.
Obviously, harmonic-Einstein and almost harmonic-Einstein are natural generalizations of Einstein metrics.

In recent years, gradient Ricci-harmonic soliton metrics and quasi Ricci-harmonic
metrics have been extensively studied by many mathematicians.
For instance, Yang and Shen \cite{Yang2012} gave a volume growth estimate for complete non-compact domain manifolds of shrinking Ricci-harmonic solitons. Tadano \cite{Tadano} gave some gap theorems for Ricci-harmonic solitons with compact domain manifolds by showing some necessary and sufficient conditions for the solitons to be harmonic-Einstein. Zhu \cite{Zhu2018} prove that when $\alpha> 0$ and the sectional curvature of $N$ is bounded from above by $\frac{\alpha}{m}$, any shrinking or steady Ricci-harmonic soliton must be a Ricci soliton. Wang \cite{Wanglf2017} studied $f$-non-parabolic ends and connectivity respectively for quasi Ricci-harmonic metrics, for more details see \cite{Guo2015,Ma2013,Wu2017,Zeng2018}.

In \cite{Barros,Huang}, the authors got some rigid properties for generalized $\tau$-quasi-Einstein metrics by establishing some integral formulas. In the first part of this paper, we derive an integral formula with conditions on the curvature to produce sufficient conditions to conclude that the compact generalized $\tau$-quasi Ricci-harmonic metric is harmonic-Einstein.

\begin{thm}\label{thm1}
Let $(M^m,g)$ be an $m$-dimensional compact manifold with a metric $g$ satisfying generalized $\tau$-quasi Ricci-harmonic metric equation \eqref{general1} with $m\geq3$. We also assume that
\begin{equation}\label{general6}
\int\limits_{M^m}\langle \nabla R_{\phi}, \nabla f\rangle\,e^{-\frac{f}{\tau}}d\upsilon\leq 0,
\end{equation}
where $d\upsilon$ is the volume form of $(M^m,g)$. Then, $(M^m,g)$ is harmonic-Einstein.
\end{thm}

In \cite{Huang13,Wanglf2018}, the authors proved that a given condition forces $(\tau, \rho)$-quasi Einstein metrics to be a standard Einstein metric. In the second part of this paper, we find rigid properties for $(\tau, \rho)$-quasi Ricci-harmonic metrics  under various conditions on $R_{\phi}$, $\rho$ and $\tau$.

\begin{thm}\label{thm2}
Let $(M^m,g)$ be an $m$-dimensional compact manifold with a metric $g$ satisfying $(\tau, \rho)$-quasi Ricci-harmonic metric equation \eqref{general2} with $m\geq3$. Then, we have
\begin{enumerate}
  \item If $\rho\geq\frac{1}{m}$, then $(M^m,g)$ is harmonic-Einstein;
  \item If $\frac{1}{2(m-1)}\leq\rho<\frac{1}{m}$ and $\tau\geq1$, then either
  $$\mu<0, ~~~~~~~~R_{\phi}\leq\frac{m(m-1)\mu}{\tau+(m-1)(1-m\rho)},$$ or $(M^m,g)$ is harmonic-Einstein;
  \item If $\rho\leq\frac{1}{2(m-1)}$ and $\tau\geq1$, then either $$\mu>0, ~~~~~~~~R_{\phi}\geq\frac{m(m-1)\mu}{\tau+(m-1)(1-m\rho)},$$ or $(M^m,g)$ is harmonic-Einstein;
  \item If $\rho=\frac{1}{2(m-1)}$ and $\tau\geq1$, then $(M^m,g)$ is harmonic-Einstein.
\end{enumerate}
\end{thm}

In the third part of this paper, we shall extend gap theorems for compact gradient Ricci solitons \cite{Fer2012},
for compact Ricci-harmonic solitons \cite{Tadano} and for compact $\tau$-quasi-Einstein metrics \cite{Wanglinf2018} to the case of compact $\tau$-quasi Ricci-harmonic metrics.

\begin{thm}\label{thm3}
Let $(M^m,g)$ be an $m$-dimensional compact manifold with a metric $g$ satisfying $\tau$-quasi Ricci-harmonic metric equation \eqref{general3} with $\lambda>0$ and $\tau>1$. Then
\begin{equation}\label{general7}
(R_{\phi})_{max}-m\lambda\leq(\tau-1)\left(\frac{2}{m(m+\tau-1)}+\frac{1}{\tau}\right)\frac{1}{V_{\tau,f}}\int\limits_{M^m}|\nabla f|^{2}\,e^{-\frac{f}{\tau}}d\upsilon,
\end{equation}
if and only if $(M^m,g)$ is harmonic-Einstein, where $(R_{\phi})_{max}$ denotes the maximal value of $R_{\phi}$ on $M$ and $V_{\tau,f}=\int\limits_{M^m}\,e^{-\frac{f}{\tau}}d\upsilon$.
\end{thm}

We point out that if $(M^m,g)$ be an $m$-dimensional compact manifold with a metric $g$ satisfying generalized $\tau$-quasi Ricci-harmonic metric equation \eqref{general1} with $\tau>1$, $\lambda>0$ and $\Delta \lambda\leq0(\geq0)$, then Theorem \ref{thm3} remains valid by Hopf lemma. For the almost Ricci-harmonic soliton case a result corresponding to the referred theorem has been obtained in \cite{Abolarinwa19}.

\begin{thm}\label{thm4}
Let $(M^m,g)$ be an $m$-dimensional compact manifold with a metric $g$ satisfying $\tau$-quasi Ricci-harmonic metric equation \eqref{general3} with $\lambda>0$ and $\tau>64m$. If
\begin{equation}\label{general8}
Ric_{\phi}\geq(1-\delta)\lambda g
\end{equation}
holds for any $\delta$ satisfying
\begin{equation}\label{general9}
0\leq\delta<\min\left\{\frac{1}{6}, \frac{\int\limits_{M^m}|\nabla f|^{2}\,e^{-\frac{f}{\tau}}d\upsilon}{3m\tau\lambda\int\limits_{M^m}\,e^{-\frac{f}{\tau}}d\upsilon}\right\}.
\end{equation}
Then, $(M^m,g)$ is harmonic-Einstein.
\end{thm}

This paper is organized as follows. In Section 2, we shall give some lemmas playing important
roles in proving Theorems \ref{thm1}-\ref{thm4}. In Section 3, a proof of Theorem \ref{thm1} shall be given. In Section 4, a proof of Theorem \ref{thm2} shall be given. Theorems \ref{thm3}-\ref{thm4} are proved in Sections 5, respectively.

\section{Preliminaries}

In this section, we first derive some basic formulas which will be
used later. Moreover, $\rho,\mu$ always stand for two
real constants in this paper unless a special explanation.
Recall that for smooth function $f$, the following operator:
\begin{align*}
\Delta_{\tau,f}(u)=&e^{\frac{f}{\tau}}div(e^{-\frac{f}{\tau}}\nabla
u)=~\Delta u-\frac{1}{\tau}\langle\nabla f,\nabla u\rangle
\end{align*}
is self-adjoint with respect to the $L^2$ inner product under the
measure $e^{-\frac{f}{\tau}}d\upsilon$ (see Lemma 3.1 in \cite{Huang14}), where $d\upsilon$ is the volume form of $(M^m,g)$. That is, $\forall \ u,v\in
C^\infty(M^m)$,
$$\int\limits_{M^m}u\Delta_{\tau,f}(v)\,e^{-\frac{f}{\tau}}d\upsilon=
-\int\limits_{M^m}\langle\nabla u,\nabla v\rangle
\,e^{-\frac{f}{\tau}}d\upsilon=\int\limits_{M^m}v\Delta_{\tau,f}(u)\,e^{-\frac{f}{\tau}}d\upsilon.$$

\begin{lem}\label{lem2-1}
Let $g$ be a generalized $\tau$-quasi Ricci-harmonic metric defined in Definition \ref{def1}. Then one can get
\begin{align}
&\Delta_{\tau,f}(f)=\Delta f-\frac{1}{\tau}|\nabla f|^2=m\lambda-R_{\phi},\label{lem2-1-1}\\
&\frac{1}{2}\nabla R_{\phi}=(m-1)\nabla\lambda+\frac{\tau-1}{\tau}Ric_{\phi}(\nabla f)+\frac{1}{\tau}[R_{\phi}-(m-1)\lambda]\nabla f,\label{lem2-1-2}\\
&\frac{1}{2}\Delta R_{\phi}=(m-1)\Delta\lambda-\frac{2(m-1)}{\tau}\langle\nabla
f,\nabla\lambda\rangle +\frac{\tau+2}{2\tau}\langle\nabla f,\nabla R_{\phi}\rangle\label{lem2-1-3}\\
&\qquad\qquad-\frac{m+\tau-1}{m\tau}(R_{\phi}-m\lambda)\left(R_{\phi}-\frac{m(m-1)}{m+\tau-1}\lambda\right)\nonumber\\
&\qquad\qquad-\frac{\tau-1}{\tau}\bigg|Ric_{\phi}-\frac{R_{\phi}}{m}g\bigg|^{2}-\frac{\tau-1}{\tau}\alpha(\tau_{g}\phi)^{2}.\nonumber
\end{align}
\end{lem}

\proof The equation \eqref{lem2-1-1} is direct consequence of the first equation in \eqref{general1}.

We use the second contracted Bianchi identity
$$\nabla R=2div Ric$$
as well as the fact that
\begin{equation}\label{lem2-1pf1}
\nabla|\nabla \phi|^{2}=2\nabla_{\nabla \phi}\nabla \phi
\end{equation}
and
$$ div(\nabla \phi \otimes \nabla \phi)=\tau_{g}\phi \nabla\phi+\nabla_{\nabla \phi}\nabla \phi$$
to deduce
\begin{equation}\label{lem2-1pf2}\aligned
\nabla R_{\phi}=&\nabla R-\alpha\nabla|\nabla \phi|^{2}\\
=&2div Ric-2\alpha\nabla_{\nabla \phi}\nabla \phi\\
=&2div Ric-2\alpha(div(\nabla \phi \otimes \nabla \phi)-\tau_{g}\phi \nabla\phi)\\
=&2div Ric_{\phi}+2\alpha\tau_{g}\phi \nabla\phi.
\endaligned\end{equation}
Using the first equation in \eqref{general1}, the Ricci identity
$$div \nabla^{2} f = Ric (\nabla f) + \nabla \Delta f$$
and $$ div(\nabla f \otimes \nabla f)=\Delta f \nabla f +\nabla_{\nabla f}\nabla f,$$
we have
\begin{equation}\label{lem2-1pf3}\aligned
\nabla R_{\phi}=&2div \left(\lambda g+\frac{1}{\tau}df\otimes df-\nabla^{2} f \right)+2\alpha\tau_{g}\phi \nabla\phi\\
=&2\nabla \lambda+\frac{2}{\tau}(\Delta f \nabla f +\nabla_{\nabla f}\nabla f)+2\alpha\tau_{g}\phi \nabla\phi\\
&-2(Ric (\nabla f) + \nabla \Delta f)\\
=&2\nabla \lambda+\frac{2}{\tau}\Delta f \nabla f+\frac{2}{\tau}\nabla_{\nabla f}\nabla f-2\nabla \Delta f-2Ric_{\phi}(\nabla f),
\endaligned\end{equation}
here we have used the second equation in \eqref{general1} in the last equality.
Using the equation \eqref{lem2-1-1} yields
\begin{equation}\label{lem2-1pf4}
\Delta f=-R_{\phi}+\lambda m+\frac{1}{\tau}|\nabla f|^{2}.
\end{equation}
Substituting \eqref{lem2-1pf4} and remembering that $\nabla|\nabla f|^{2}=2\nabla_{\nabla f}\nabla f$ we use \eqref{lem2-1pf3} to write
\begin{equation}\label{lem2-1pf5}\aligned
\nabla R_{\phi}=&2\nabla \lambda+\frac{2}{\tau}\Delta f \nabla f+\frac{2}{\tau}\nabla_{\nabla f}\nabla f-2Ric_{\phi}(\nabla f)\\
&-2\nabla (-R_{\phi}+\lambda m+\frac{1}{\tau}|\nabla f|^{2})\\
=&-2(m-1)\nabla \lambda-2Ric_{\phi}(\nabla f)+2\nabla R_{\phi}\\
&+\frac{2}{\tau}\Delta f \nabla f-\frac{2}{\tau}\nabla_{\nabla f}\nabla f.
\endaligned\end{equation}
We now use \eqref{lem2-1pf5} to write
\begin{equation}\label{lem2-1pf6}
\nabla R_{\phi}=2(m-1)\nabla \lambda+\frac{2}{\tau}\nabla_{\nabla f}\nabla f-\frac{2}{\tau}\Delta f \nabla f+2Ric_{\phi}(\nabla f).
\end{equation}
From the first equation of \eqref{general1}, we have
\begin{equation}\label{lem2-1pf7}
\nabla_{\nabla f}\nabla f=\frac{1}{\tau}|\nabla f|^{2}\nabla f+\lambda\nabla f-Ric_{\phi}(\nabla f).
\end{equation}
Insertting \eqref{lem2-1pf4} and \eqref{lem2-1pf7} into \eqref{lem2-1pf6} leads to \eqref{lem2-1-2}, which finishes the second statement of the lemma.

Initially by using \eqref{lem2-1-2} to compute the divergence of $\nabla R_{\phi}$ we obtain
\begin{equation}\label{lem2-1pf8}\aligned
\frac{1}{2}\Delta R_{\phi}=&(m-1)\Delta \lambda+\frac{\tau-1}{\tau}div(Ric_{\phi}(\nabla f))\\
&+\frac{1}{\tau}div((R_{\phi}-(m-1)\lambda)\nabla f).
\endaligned\end{equation}
By \eqref{lem2-1pf2} we have
\begin{equation}\label{lem2-1pf9}\aligned
div(Ric_{\phi}(\nabla f))
=&div(Ric_{\phi})(\nabla f)+trace(Ric_{\phi}\circ\nabla^{2}f)\\
=&trace\left(Ric_{\phi}\circ\left(\frac{1}{\tau}df\otimes df+\lambda g-Ric_{\phi}\right)\right)\\
&+\left(\frac{1}{2}\nabla R_{\phi}-\alpha\tau_{g}\phi \nabla\phi\right)(\nabla f)\\
=&\frac{1}{2}\langle \nabla R_{\phi}, \nabla f\rangle-\alpha(\tau_{g}\phi)^{2}+\frac{1}{\tau}Ric_{\phi}(\nabla f,\nabla f)\\
&+trace\left(Ric_{\phi}\circ\left(\lambda g-Ric_{\phi}\right)\right).
\endaligned\end{equation}
By \eqref{lem2-1-2} we have
\begin{equation}\label{lem2-1pf10}\aligned
\frac{1}{\tau}Ric_{\phi}(\nabla f,\nabla f)
=&\frac{1}{2(\tau-1)}\langle \nabla R_{\phi}, \nabla f\rangle-\frac{m-1}{\tau-1}\langle \nabla \lambda, \nabla f\rangle\\
&-\frac{1}{\tau(\tau-1)}(R_{\phi}-(m-1)\lambda)|\nabla f|^{2}.
\endaligned\end{equation}
Plugging \eqref{lem2-1pf10} into \eqref{lem2-1pf9} leads to
\begin{equation}\label{lem2-1pf11}\aligned
div(Ric_{\phi}(\nabla f))
=&\frac{\tau}{2(\tau-1)}\langle \nabla R_{\phi}, \nabla f\rangle-\alpha(\tau_{g}\phi)^{2}-\frac{m-1}{\tau-1}\langle \nabla \lambda, \nabla f\rangle\\
&-\frac{1}{\tau(\tau-1)}(R_{\phi}-(m-1)\lambda)|\nabla f|^{2}\\
&+trace\left(Ric_{\phi}\circ\left(\lambda g-Ric_{\phi}\right)\right).
\endaligned\end{equation}
Insertting \eqref{lem2-1pf11} into \eqref{lem2-1pf8} leads to
\begin{equation}\label{lem2-1pf12}\aligned
\frac{1}{2}\Delta R_{\phi}=&(m-1)\Delta\lambda-\frac{2(m-1)}{\tau}\langle\nabla
f,\nabla\lambda\rangle +\frac{\tau+2}{2\tau}\langle\nabla f,\nabla R_{\phi}\rangle\\
&-\frac{\tau-1}{\tau}\alpha(\tau_{g}\phi)^{2}+\frac{\tau-1}{\tau}trace\left(Ric_{\phi}\circ\left(\lambda g-Ric_{\phi}\right)\right)\\
&+\frac{1}{\tau}(R_{\phi}-(m-1)\lambda)(m\lambda-R_{\phi}).
\endaligned\end{equation}
Note that
\begin{equation}\label{lem2-1pf13}\aligned
&trace\left(Ric_{\phi}\circ\left(\lambda g-Ric_{\phi}\right)\right)=\lambda R_{\phi}-|Ric_{\phi}|^{2}\\
=&-\bigg|Ric_{\phi}-\frac{R_{\phi}}{m}g\bigg|^{2}+R_{\phi}\left(\lambda-\frac{1}{m}R_{\phi}\right).
\endaligned\end{equation}
Plugging \eqref{lem2-1pf13} into \eqref{lem2-1pf12}, we arrive at \eqref{lem2-1-3}.
\endproof

\begin{lem}\label{lem2-2}
If $g$ is a generalized $\tau$-quasi Ricci-harmonic metric defined in Definition \ref{def1} and $\lambda=F(f)$, where $F(t)$ is a smooth function, then there exists a constant $\zeta$, so that
\begin{equation}\label{lem2-2-1}
R_{\phi}+\frac{\tau-1}{\tau}|\nabla f|^{2}=\left[2(m-1)G^{'}(f)+\frac{2(m+\tau-2)}{\tau}G(f)+\zeta\right]e^{\frac{2}{\tau}f},
\end{equation}
where $$G(t)=\int F(t)e^{-\frac{2}{\tau}t}\,dt.$$
\end{lem}

\proof Using \eqref{lem2-1-2} and the first equation of \eqref{general1} we can write
\begin{equation}\label{lem2-2pf1}\aligned
\nabla R_{\phi}=&\frac{2(\tau-1)}{\tau}\left(\frac{1}{\tau}df\otimes df+\lambda g-\nabla^{2}f\right)(\nabla f)\\
&+2(m-1)\nabla \lambda+\frac{2}{\tau}R_{\phi}\nabla f-\frac{2(m-1)\lambda}{\tau}\nabla f\\
=&\frac{2(\tau-1)}{\tau}\lambda\nabla f+\frac{2(\tau-1)}{\tau^{2}}|\nabla f|^{2}\nabla f-\frac{2(\tau-1)}{\tau}\nabla^{2}f(\nabla f)\\
&+2(m-1)\nabla \lambda+\frac{2}{\tau}R_{\phi}\nabla f-\frac{2(m-1)}{\tau}\lambda\nabla f\\
=&\frac{2(\tau-m)}{\tau}\lambda\nabla f+2(m-1)\nabla \lambda+\frac{2}{\tau}R_{\phi}\nabla f\\
&+\frac{2(\tau-1)}{\tau^{2}}|\nabla f|^{2}\nabla f-\frac{\tau-1}{\tau}\nabla|\nabla f|^{2}.
\endaligned\end{equation}
By the definition of $G$ and the fact that $\nabla \lambda=F^{'}(f)\nabla f$, \eqref{lem2-2pf1} can be rewritten as
$$\nabla\left[\left(R_{\phi}+\frac{\tau-1}{\tau}|\nabla f|^{2}\right)e^{-\frac{2}{\tau}f}-2(m-1)G^{'}(f)-\frac{2(m+\tau-2)}{\tau}G(f)\right]=0$$
Therefore,
$$\left(R_{\phi}+\frac{\tau-1}{\tau}|\nabla f|^{2}\right)e^{-\frac{2}{\tau}f}-2(m-1)G^{'}(f)-\frac{2(m+\tau-2)}{\tau}G(f)$$ is constant, and \eqref{lem2-2-1} follows.
\endproof

\begin{lem}\label{lem2-3}
Let $g$ be a $(\tau, \rho)$-quasi Ricci-harmonic metric defined in Definition \ref{def2}. Then one can get
\begin{align}
&\Delta_{\tau,f}(f)=\Delta f-\frac{1}{\tau}|\nabla f|^2=(m\rho-1)R_{\phi}+m\mu,\label{lem2-3-1}\\
&\frac{2(m-1)\rho-1}{2}\Delta R_{\phi}=\frac{4(m-1)\rho-(\tau+2)}{2\tau}\langle\nabla f,\nabla R_{\phi}\rangle\label{lem2-3-2}\\
&\qquad\qquad+\frac{(1-m\rho)R_{\phi}-m\mu}{m\tau}\big[(\tau+(m-1)(1-m\rho))R_{\phi}-m(m-1)\mu\big]\nonumber\\
&\qquad\qquad+\frac{\tau-1}{\tau}\bigg|Ric_{\phi}-\frac{R_{\phi}}{m}g\bigg|^{2}+\frac{\tau-1}{\tau}\alpha(\tau_{g}\phi)^{2}.\nonumber
\end{align}
\end{lem}

\proof Taking $\lambda=\rho R+\mu$ in the equations \eqref{lem2-1-1} and \eqref{lem2-1-3}, we obtain
the formulas \eqref{lem2-3-1} and \eqref{lem2-3-2}.
\endproof

\begin{lem}\label{lem2-4}(\cite{Wanglf}) Let $g$ be a $\tau$-quasi Ricci-harmonic metric defined in Definition \ref{def3}. Then one can get
\begin{equation}\label{lem2-4-0}
\Delta_{\tau,f}(f)=\Delta f-\frac{1}{\tau}|\nabla f|^2=m\lambda-R_{\phi},
\end{equation}

\begin{equation}\label{lem2-4-1}
\frac{1}{2}\nabla R_{\phi}=\frac{\tau-1}{\tau}Ric_{\phi}(\nabla f)+\frac{1}{\tau}[R_{\phi}-(m-1)\lambda]\nabla f
\end{equation}
and
\begin{equation}\label{lem2-4-2}\aligned
\frac{1}{2}\Delta R_{\phi}=&\frac{\tau+2}{2\tau}\langle\nabla R_{\phi}, \nabla f\rangle-\frac{\tau-1}{\tau}\alpha(\tau_{g}\phi)^{2}-\frac{\tau-1}{\tau}\left|Ric_{\phi}-\frac{1}{m}R_{\phi}g\right|^{2}\\
&-\frac{m+\tau-1}{m\tau}(R_{\phi}-m\lambda)\left(R_{\phi}-\frac{m(m-1)}{m+\tau-1}\lambda\right).
\endaligned\end{equation}

Moreover, there exists a constant $\varrho$ such that

\begin{equation}\label{lem2-4-3}
R_{\phi}+\frac{\tau-1}{\tau}|\nabla f|^{2}-(m-\tau)\lambda=\varrho e^{\frac{2}{\tau}f}
\end{equation}
and
\begin{equation}\label{lem2-4-4}
\Delta f + |\nabla f|^{2}-\tau\lambda+\varrho e^{\frac{2}{\tau}f}=0.
\end{equation}
\end{lem}

\proof Since $\lambda$ is constant, then from Lemmas \ref{lem2-1} and \ref{lem2-2} we conclude the proof of the lemma.
\endproof

\section{Proof of Theorem \ref{thm1}}

In this section, we prove Theorem \ref{thm1}. Firstly, we introduce an integral formula for a compact generalized $\tau$-quasi Ricci-harmonic metric.

\begin{lem}\label{lem3-1}
Let $(M^m,g)$ be an $m$-dimensional compact manifold with a metric $g$ satisfying generalized $\tau$-quasi Ricci-harmonic metric equation \eqref{general1} with $m\geq3$. Then, we have
\begin{equation}\label{lem3-1-1}\aligned
\int\limits_{M^m}\bigg|Ric_{\phi}-\frac{R_{\phi}}{m}g\bigg|^{2}\,e^{-\frac{f}{\tau}}d\upsilon&+\alpha\int\limits_{M^m}(\tau_{g}\phi)^{2}\,e^{-\frac{f}{\tau}}d\upsilon\\
&=\frac{m-2}{2m}\int\limits_{M^m}\langle\nabla R_{\phi},\nabla f\rangle\,e^{-\frac{f}{\tau}}d\upsilon.
\endaligned\end{equation}
\end{lem}

\begin{rem}\label{rem3-1}
If $\phi: (M, g) \to (\mathbb{R}, dr^{2})$ is a constant function in \eqref{lem3-1-1}, we obtain
the formula (3.1) in \cite{Huang}.
\end{rem}

\proof  We have from \eqref{lem2-1-1}
$$\Delta_{\tau,f}(f)=m\lambda-R_{\phi}.$$
In particular, \eqref{lem2-1-3} can be written as
$$\aligned
&\Delta_{\tau,f}\left(\frac{1}{2}R_{\phi}-(m-1)\lambda\right)\\
=&\frac{1}{2}\left(\Delta R_{\phi}-\frac{1}{\tau}\langle\nabla f,\nabla R_{\phi}\rangle \right)-(m-1)\left(\Delta \lambda-\frac{1}{\tau}\langle\nabla f,\nabla \lambda\rangle \right)\\
=&-\frac{m-1}{\tau}\langle\nabla f,\nabla \lambda\rangle +\frac{\tau+1}{2\tau}\langle\nabla f,\nabla R_{\phi}\rangle-\frac{\tau-1}{\tau}\alpha(\tau_{g}\phi)^{2}\\
&-\frac{\tau-1}{\tau}\bigg|Ric_{\phi}-\frac{R_{\phi}}{m}g\bigg|^{2}
-\frac{(R_{\phi}-m\lambda)[(m+\tau-1)R_{\phi}-m(m-1)\lambda]}{m\tau}.\endaligned$$
Integrating the above formula yields
$$\aligned
0=&\int\limits_{M^m}\biggl\{-\frac{m-1}{\tau}\langle\nabla f,\nabla \lambda\rangle +\frac{\tau+1}{2\tau}\langle\nabla f,\nabla R_{\phi}\rangle
-\frac{\tau-1}{\tau}\bigg|Ric_{\phi}-\frac{R_{\phi}}{m}g\bigg|^{2}\\
&-\frac{\tau-1}{\tau}\alpha(\tau_{g}\phi)^{2}-\frac{(R_{\phi}-m\lambda)[(m+\tau-1)R_{\phi}-m(m-1)\lambda]}{m\tau}\biggr\}\,e^{-\frac{f}{\tau}}d\upsilon\\
=&\int\limits_{M^m}\biggl\{\frac{m-1}{\tau}\lambda\Delta_{\tau,f}(f)
-\frac{\tau+1}{2\tau}R_{\phi}\,\Delta_{\tau,f}(f)
-\frac{\tau-1}{\tau}\bigg|Ric_{\phi}-\frac{R_{\phi}}{m}g\bigg|^{2}\\
&-\frac{\tau-1}{\tau}\alpha(\tau_{g}\phi)^{2}-\frac{(R_{\phi}-m\lambda)[(m+\tau-1)R_{\phi}-m(m-1)\lambda]}{m\tau}\biggr\}\,e^{-\frac{f}{\tau}}d\upsilon\\
=&\int\limits_{M^m}\biggl\{-\frac{(m-2)(\tau-1)}{2m\tau}R_{\phi}(m\lambda-R_{\phi})-\frac{\tau-1}{\tau}\alpha(\tau_{g}\phi)^{2}\\
&\,\,\,\,\,\,\,\,\,\,\,\,\,\,\,\,\,\,\,\,\,\,\,\,\,\,\,\,\,\,\,\,\,\,\,\,\,\,\,\,\,\,\,\,\,\,\,\,\,\,\,\,\,\,\,\,\,\,\,\,\,\,\,\,\,\,\,\,\,\,\,\,\,\,\,\,\,\,\,\,\,\,\,\,\,\,\,\,\,\,\,\,\,\,\,\,
-\frac{\tau-1}{\tau}\bigg|Ric_{\phi}-\frac{R_{\phi}}{m}g\bigg|^{2}\biggr\}\,e^{-\frac{f}{\tau}}d\upsilon\\
=&\int\limits_{M^m}\biggl\{-\frac{(m-2)(\tau-1)}{2m\tau}R_{\phi}\Delta_{\tau,f}(f)-\frac{\tau-1}{\tau}\alpha(\tau_{g}\phi)^{2}\\
&\,\,\,\,\,\,\,\,\,\,\,\,\,\,\,\,\,\,\,\,\,\,\,\,\,\,\,\,\,\,\,\,\,\,\,\,\,\,\,\,\,\,\,\,\,\,\,\,\,\,\,\,\,\,\,\,\,\,\,\,\,\,\,\,\,\,\,\,\,\,\,\,\,\,\,\,\,\,\,\,\,\,\,\,\,\,\,\,\,\,\,\,\,\,\,\,
-\frac{\tau-1}{\tau}\bigg|Ric_{\phi}-\frac{R_{\phi}}{m}g\bigg|^{2}\biggr\}\,e^{-\frac{f}{\tau}}d\upsilon\\
=&\int\limits_{M^m}\biggl\{\frac{(m-2)(\tau-1)}{2m\tau}\langle\nabla R_{\phi},\nabla f\rangle-\frac{\tau-1}{\tau}\alpha(\tau_{g}\phi)^{2}\\
&\,\,\,\,\,\,\,\,\,\,\,\,\,\,\,\,\,\,\,\,\,\,\,\,\,\,\,\,\,\,\,\,\,\,\,\,\,\,\,\,\,\,\,\,\,\,\,\,\,\,\,\,\,\,\,\,\,\,\,\,\,\,\,\,\,\,\,\,\,\,\,\,\,\,\,\,\,\,\,\,\,\,\,\,\,\,\,\,\,\,\,\,\,\,\,\,
-\frac{\tau-1}{\tau}\bigg|Ric_{\phi}-\frac{R_{\phi}}{m}g\bigg|^{2}\biggr\}\,e^{-\frac{f}{\tau}}d\upsilon.
\endaligned$$
We complete the proof of Lemma \ref{lem3-1}.\endproof

Now we prove Theorem \ref{thm1}.\\
\textit{Proof of Theorem}\ref{thm1}.
Since $m\geq 3$ and the condition $$\int\limits_{M^m}\langle \nabla R_{\phi}, \nabla f\rangle\,e^{-\frac{f}{\tau}}d\upsilon\leq 0,$$ then from Lemma \ref{lem3-1} we obtain
$$\int\limits_{M^m}\bigg|Ric_{\phi}-\frac{R_{\phi}}{m}g\bigg|^{2}\,e^{-\frac{f}{\tau}}d\upsilon+\alpha\int\limits_{M^m}(\tau_{g}\phi)^{2}\,e^{-\frac{f}{\tau}}d\upsilon=0,$$
which gives $Ric_{\phi}-\frac{R_{\phi}}{m}g=0$ and $\tau_{g}\phi=0$, completing the proof of the theorem.

\section{Proof of Theorem \ref{thm2}}

In this section, we study rigid properties of a special generalized $\tau$-quasi Ricci-harmonic metric with $\lambda=\rho R_{\phi}+\mu$, where
$\rho,\mu$ are two real constants. First, we will use the following important lemma.

\begin{lem}\label{lem4-1}
Let $(M^m,g)$ be an $m$-dimensional compact manifold with a metric $g$ satisfying  $(\tau, \rho)$-quasi Ricci-harmonic metric equation \eqref{general2}. If the generalized scalar curvature $R_{\phi}$ is constant, then $(M^m,g)$ is harmonic-Einstein.
\end{lem}

\proof Integrating \eqref{lem2-3-1} with respect to $e^{-\frac{f}{\tau}}d\upsilon$ gives
$$\aligned
((m\rho-1)R_{\phi}+m\mu)\int\limits_{M^m}\,e^{-\frac{f}{\tau}}d\upsilon=&\int\limits_{M^m}\left(\Delta f-\frac{1}{\tau}|\nabla f|^2\right)\,e^{-\frac{f}{\tau}}d\upsilon\\
=&\int\limits_{M^m}\Delta_{\tau,f}(f)\,e^{-\frac{f}{\tau}}d\upsilon\\
=&0.
\endaligned$$
Hence $(m\rho-1)R_{\phi}+m\mu = 0$ and then
\begin{equation}\label{lem4-1-pf1}
\Delta f-\frac{1}{\tau}|\nabla f|^2=0.
\end{equation}
Integrating \eqref{lem4-1-pf1} with respect to $d\upsilon$ gives
$$\frac{1}{\tau}\int\limits_{M^m}|\nabla f|^2\,d\upsilon=\int\limits_{M^m}\Delta f\,d\upsilon=0,$$
which gives $f$ is a constant and $(M^m,g)$ is harmonic-Einstein.
\endproof

Now we prove Theorem \ref{thm2}.\\
\textit{Proof of Theorem}\ref{thm2}. \textbf{(1)} If $\rho=\frac{1}{m}$, then \eqref{lem2-3-1} shows
$$\Delta_{\tau,f}(f)=\Delta f-\frac{1}{\tau}|\nabla f|^2=m\mu,$$
which gives that $\mu = 0$ and $f$ is constant and $(M^m,g)$ is harmonic-Einstein.

Now we consider the case $\rho>\frac{1}{m}$. Integrating \eqref{lem2-3-1} on $M$ with respect to the measure $e^{-\frac{f}{\tau}}d\upsilon$ leads to
\begin{equation}\label{thm2-pf1}
m\mu\int\limits_{M^m}\,e^{-\frac{f}{\tau}}d\upsilon=(1-m\rho)\int\limits_{M^m}R_{\phi}\,e^{-\frac{f}{\tau}}d\upsilon.
\end{equation}
Integrating by parts and using the first equation in \eqref{general2}, we have
$$\aligned
&\int\limits_{M^m}Ric_{\phi}(\nabla f, \nabla f)\,e^{-\frac{f}{\tau}}d\upsilon\\
=&-\tau\int\limits_{M^m}Ric_{\phi}(\nabla f, \nabla e^{-\frac{f}{\tau}})\,d\upsilon\\
=&\tau\int\limits_{M^m}\{(div Ric_{\phi})(\nabla f)+trace(Ric_{\phi}\circ\nabla^{2}f)\}\,e^{-\frac{f}{\tau}}d\upsilon\\
=&\tau\int\limits_{M^m}\biggl\{\big\langle \frac{1}{2}\nabla R_{\phi}, \nabla f\big\rangle-\alpha(\tau_{g}\phi)^{2}\\
&\,\,\,\,\,\,\,\,\,\,\,\,\,
+trace\left(Ric_{\phi}\circ\big(\rho R_{\phi}g+\mu g-Ric_{\phi}+\frac{1}{\tau}df\otimes df\big)\right)\biggr\}\,e^{-\frac{f}{\tau}}d\upsilon\\
=&\tau\int\limits_{M^m}\biggl\{\big\langle \frac{1}{2}\nabla R_{\phi}, \nabla f\big\rangle-\alpha(\tau_{g}\phi)^{2}\\
&\,\,\,\,\,\,\,\,\,\,\,\,\,
+\rho R_{\phi}^{2}+\mu R_{\phi}+\frac{1}{\tau}Ric_{\phi}(\nabla f, \nabla f)-|Ric_{\phi}|^{2}\biggr\}\,e^{-\frac{f}{\tau}}d\upsilon.
\endaligned$$
Hence
$$\aligned
0=&\int\limits_{M^m}\biggl\{\big\langle \frac{1}{2}\nabla R_{\phi}, \nabla f\big\rangle-\alpha(\tau_{g}\phi)^{2}+\rho R_{\phi}^{2}+\mu R_{\phi}-|Ric_{\phi}|^{2}\biggr\}\,e^{-\frac{f}{\tau}}d\upsilon\\
=&\int\limits_{M^m}\biggl\{-\frac{1}{2}R_{\phi}\Delta_{\tau,f}(f)-\alpha(\tau_{g}\phi)^{2}+\rho R_{\phi}^{2}+\mu R_{\phi}-|Ric_{\phi}|^{2}\biggr\}\,e^{-\frac{f}{\tau}}d\upsilon\\
=&\int\limits_{M^m}\biggl\{-\frac{1}{2}R_{\phi}((m\rho-1)R_{\phi}+m\mu)-\alpha(\tau_{g}\phi)^{2}+\rho R_{\phi}^{2}+\mu R_{\phi}-|Ric_{\phi}|^{2}\biggr\}\,e^{-\frac{f}{\tau}}d\upsilon.
\endaligned$$
Note that
$$|Ric_{\phi}|^{2}\geq\frac{1}{m}R_{\phi}^{2}.$$
We then get
\begin{equation}\label{thm2-pf2}\aligned
0\leq&\int\limits_{M^m}\biggl\{-\frac{1}{2}R_{\phi}((m\rho-1)R_{\phi}+m\mu)-\alpha(\tau_{g}\phi)^{2}+\rho R_{\phi}^{2}+\mu R_{\phi}-\frac{1}{m}R_{\phi}^{2}\biggr\}\,e^{-\frac{f}{\tau}}d\upsilon\\
=&\frac{m-2}{2}\int\limits_{M^m}\biggl(\frac{1-m\rho}{m}R_{\phi}^{2}-\mu R_{\phi}\biggr)\,e^{-\frac{f}{\tau}}d\upsilon-\alpha\int\limits_{M^m}(\tau_{g}\phi)^{2}\,e^{-\frac{f}{\tau}}d\upsilon\\
=&\frac{(m-2)(1-m\rho)}{2m}\left(\int\limits_{M^m}R_{\phi}^{2}\,e^{-\frac{f}{\tau}}d\upsilon-\frac{\left(\int\limits_{M^m}R_{\phi}\,e^{-\frac{f}{\tau}}d\upsilon\right)^{2}}{\int\limits_{M^m}\,e^{-\frac{f}{\tau}}d\upsilon} \right)-\alpha\int\limits_{M^m}(\tau_{g}\phi)^{2}\,e^{-\frac{f}{\tau}}d\upsilon,
\endaligned\end{equation}
where we have used \eqref{thm2-pf1}. Since $1-m\rho<0$, then \eqref{thm2-pf2} implies
$$\aligned
\left(\int\limits_{M^m}R_{\phi}\,e^{-\frac{f}{\tau}}d\upsilon\right)^{2}\geq&\int\limits_{M^m}\,e^{-\frac{f}{\tau}}d\upsilon\,\int\limits_{M^m}R_{\phi}^{2}\,e^{-\frac{f}{\tau}}d\upsilon\\
&-\frac{2m\alpha}{(m-2)(1-m\rho)}\int\limits_{M^m}\,e^{-\frac{f}{\tau}}d\upsilon\,\int\limits_{M^m}(\tau_{g}\phi)^{2}\,e^{-\frac{f}{\tau}}d\upsilon,
\endaligned$$
which gives
\begin{equation}\label{thm2-pf3}
\left(\int\limits_{M^m}R_{\phi}\,e^{-\frac{f}{\tau}}d\upsilon\right)^{2}\geq\int\limits_{M^m}\,e^{-\frac{f}{\tau}}d\upsilon\,\int\limits_{M^m}R_{\phi}^{2}\,e^{-\frac{f}{\tau}}d\upsilon.
\end{equation}
We use \eqref{thm2-pf3} and the Cauchy-Schwartz inequality
$$\left(\int\limits_{M^m}R_{\phi}\,e^{-\frac{f}{\tau}}d\upsilon\right)^{2}\leq\int\limits_{M^m}\,e^{-\frac{f}{\tau}}d\upsilon\,\int\limits_{M^m}R_{\phi}^{2}\,e^{-\frac{f}{\tau}}d\upsilon$$
to obtain the relation
$$\left(\int\limits_{M^m}R_{\phi}\,e^{-\frac{f}{\tau}}d\upsilon\right)^{2}=\int\limits_{M^m}\,e^{-\frac{f}{\tau}}d\upsilon\,\int\limits_{M^m}R_{\phi}^{2}\,e^{-\frac{f}{\tau}}d\upsilon.$$
We then use \eqref{thm2-pf1} to obtain
$$\int\limits_{M^m}\left(R_{\phi}^{2}-\left(\frac{m\mu}{1-m\rho}\right)^{2}\right)\,e^{-\frac{f}{\tau}}d\upsilon=0,$$
which gives $R_{\phi}$ is constant and $(M^m,g)$ is harmonic-Einstein.

\textbf{(2)} Applying \eqref{lem2-3-2} to a point of maximal of $R_{\phi}$, we have
\begin{equation}\label{thm2-pf5}\aligned
&\frac{(1-m\rho)(R_{\phi})_{max}-m\mu}{m\tau}\big[(\tau+(m-1)(1-m\rho))(R_{\phi})_{max}-m(m-1)\mu\big]\\
=&\frac{2(m-1)\rho-1}{2}\Delta (R_{\phi})_{max}-\frac{\tau-1}{\tau}\bigg|Ric_{\phi}-\frac{(R_{\phi})_{max}}{m}g\bigg|^{2}-\frac{\tau-1}{\tau}\alpha(\tau_{g}\phi)^{2}\\
\leq&0.
\endaligned\end{equation}
If $\mu>0$, \eqref{thm2-pf5} implies
\begin{equation}\label{thm2-pf6}
\frac{m(m-1)\mu}{\tau+(m-1)(1-m\rho)}\leq(R_{\phi})_{max}\leq\frac{m\mu}{1-m\rho}.
\end{equation}
If $\mu\leq0$, \eqref{thm2-pf5} gives
$$\frac{m\mu}{1-m\rho}\leq(R_{\phi})_{max}\leq\frac{m(m-1)\mu}{\tau+(m-1)(1-m\rho)}.$$
Then if $R_{\phi}$ is not a constant, from \eqref{thm2-pf1} we have
$$m\mu\int\limits_{M^m}\,e^{-\frac{f}{\tau}}d\upsilon=(1-m\rho)\int\limits_{M^m}R_{\phi}\,e^{-\frac{f}{\tau}}d\upsilon<(1-m\rho)(R_{\phi})_{max}\int\limits_{M^m}\,e^{-\frac{f}{\tau}}d\upsilon,$$
which shows
$$(R_{\phi})_{max}>\frac{m\mu}{1-m\rho}.$$
This contradicts with \eqref{thm2-pf6}. Hence, $R_{\phi}$ is a constant and $(M^m,g)$ is harmonic-Einstein.

\textbf{(3)} Applying \eqref{lem2-3-2} to a point of minimum of $R_{\phi}$, we have
\begin{equation}\label{thm2-pf7}\aligned
&\frac{(1-m\rho)(R_{\phi})_{min}-m\mu}{m\tau}\big[(\tau+(m-1)(1-m\rho))(R_{\phi})_{min}-m(m-1)\mu\big]\\
=&\frac{2(m-1)\rho-1}{2}\Delta (R_{\phi})_{min}-\frac{\tau-1}{\tau}\bigg|Ric_{\phi}-\frac{(R_{\phi})_{min}}{m}g\bigg|^{2}-\frac{\tau-1}{\tau}\alpha(\tau_{g}\phi)^{2}\\
\leq&0.
\endaligned\end{equation}
If $\mu>0$, \eqref{thm2-pf7} implies
$$\frac{m(m-1)\mu}{\tau+(m-1)(1-m\rho)}\leq(R_{\phi})_{min}\leq\frac{m\mu}{1-m\rho}.$$
If $\mu\leq0$, \eqref{thm2-pf7} gives
\begin{equation}\label{thm2-pf8}
\frac{m\mu}{1-m\rho}\leq(R_{\phi})_{min}\leq\frac{m(m-1)\mu}{\tau+(m-1)(1-m\rho)}.\end{equation}
Then if $R_{\phi}$ is not a constant, from \eqref{thm2-pf1} we have
$$m\mu\int\limits_{M^m}\,e^{-\frac{f}{\tau}}d\upsilon=(1-m\rho)\int\limits_{M^m}R_{\phi}\,e^{-\frac{f}{\tau}}d\upsilon>(1-m\rho)(R_{\phi})_{min}\int\limits_{M^m}\,e^{-\frac{f}{\tau}}d\upsilon,$$
which shows
$$(R_{\phi})_{min}<\frac{m\mu}{1-m\rho}.$$
This contradicts with \eqref{thm2-pf8}. Hence, $R_{\phi}$ is a constant and $(M^m,g)$ is harmonic-Einstein.

\textbf{(4)} This result follows from \textbf{(2)} and \textbf{(3)}.

We complete the proof of the theorem.

\section{Proof of Theorems \ref{thm3} and \ref{thm4}}

In this section,
using estimates for the generalized Ricci curvature and the generalized scalar curvature, we shall give
two gap theorems for compact $\tau$-quasi Ricci-harmonic metrics by showing
some necessary and sufficient conditions for the metrics to be harmonic-Einstein. The following lemma plays crucial roles in this section:

\begin{lem}\label{lem5-1}
Let $(M^m,g)$ be an $m$-dimensional compact manifold with a metric $g$ satisfying $\tau$-quasi Ricci-harmonic metric equation \eqref{general3} with $\lambda>0$ and $\tau>1$. Then
\begin{equation}\label{lem5-1-1}\aligned
&\int\limits_{M^m}(\Delta_{\tau,f}(f))^{2}\,e^{-\frac{f}{\tau}}d\upsilon\\
=&\int\limits_{M^m}\left(\frac{2(\tau-1)}{\tau}Ric_{\phi}(\nabla f,\nabla f)+\frac{2}{\tau}[R_{\phi}-(m-1)\lambda]|\nabla f|^{2}\right)\,e^{-\frac{f}{\tau}}d\upsilon,
\endaligned\end{equation}
\begin{equation}\label{lem5-1-2}\aligned
&\int\limits_{M^m}(\Delta_{\tau,f}(f))^{2}\,e^{-\frac{f}{\tau}}d\upsilon\\
=&\int\limits_{M^m}\left(\frac{(2+m)\tau-3m}{\tau}\lambda|\nabla f|^{2}-\frac{\tau-3}{\tau}R_{\phi}|\nabla f|^{2}+\frac{2(\tau-1)}{\tau^{2}}|\nabla f|^{4}\right)\,e^{-\frac{f}{\tau}}d\upsilon
\endaligned\end{equation}
and
\begin{equation}\label{lem5-1-3}
\frac{\tau-1}{\tau}|\nabla f|^{2}\leq (R_{\phi})_{max}-R_{\phi}.
\end{equation}
\end{lem}

\proof First, by \eqref{lem2-4-0} and \eqref{lem2-4-1}, we have
\begin{equation}\label{lem5-1-pf1}\aligned
&\int\limits_{M^m}(\Delta_{\tau,f}(f))^{2}\,e^{-\frac{f}{\tau}}d\upsilon\\
=&\int\limits_{M^m}(\Delta_{\tau,f}(f))(m\lambda-R_{\phi})\,e^{-\frac{f}{\tau}}d\upsilon
=\int\limits_{M^m}\langle \nabla R_{\phi}, \nabla f\rangle\,e^{-\frac{f}{\tau}}d\upsilon\\
=&\int\limits_{M^m}\left(\frac{2(\tau-1)}{\tau}Ric_{\phi}(\nabla f,\nabla f)+\frac{2}{\tau}[R_{\phi}-(m-1)\lambda]|\nabla f|^{2}\right)\,e^{-\frac{f}{\tau}}d\upsilon,
\endaligned\end{equation}
which proves \eqref{lem5-1-1}.
Secondly, by the first equation of \eqref{general3}, we have
\begin{equation}\label{lem5-1-pf2}\aligned
&\int\limits_{M^m}(\Delta_{\tau,f}(f))^{2}\,e^{-\frac{f}{\tau}}d\upsilon\\
=&\int\limits_{M^m}\left(\frac{2(\tau-1)}{\tau}Ric_{\phi}(\nabla f,\nabla f)+\frac{2}{\tau}[R_{\phi}-(m-1)\lambda]|\nabla f|^{2}\right)\,e^{-\frac{f}{\tau}}d\upsilon\\
=&\int\limits_{M^m}\left(-\frac{\tau-1}{\tau}\langle \nabla|\nabla f|^{2}, \nabla f\rangle+\frac{2(\tau-1)}{\tau^{2}}|\nabla f|^{4}+\frac{2}{\tau}[R_{\phi}-(m-1)\lambda]|\nabla f|^{2}\right)\,e^{-\frac{f}{\tau}}d\upsilon\\
=&\int\limits_{M^m}\left(\frac{\tau-1}{\tau}|\nabla f|^{2}(m\lambda-R_{\phi})+\frac{2(\tau-1)}{\tau^{2}}|\nabla f|^{4}+\frac{2}{\tau}[R_{\phi}-(m-1)\lambda]|\nabla f|^{2}\right)\,e^{-\frac{f}{\tau}}d\upsilon\\
=&\int\limits_{M^m}\left(\frac{(2+m)\tau-3m}{\tau}\lambda|\nabla f|^{2}-\frac{\tau-3}{\tau}R_{\phi}|\nabla f|^{2}+\frac{2(\tau-1)}{\tau^{2}}|\nabla f|^{4}\right)\,e^{-\frac{f}{\tau}}d\upsilon,
\endaligned\end{equation}
where in the third equality we have used the fact that
$$\int\limits_{M^m}\langle \nabla|\nabla f|^{2}, \nabla f\rangle\,e^{-\frac{f}{\tau}}d\upsilon=-\int\limits_{M^m}|\nabla f|^{2}\Delta_{\tau,f}(f)\,e^{-\frac{f}{\tau}}d\upsilon=-\int\limits_{M^m}|\nabla f|^{2}(m\lambda-R_{\phi})\,e^{-\frac{f}{\tau}}d\upsilon.$$
Hence we arrive at \eqref{lem5-1-2}.
Finally, in order to prove \eqref{lem5-1-3}, recall from \eqref{lem2-4-3} that
$$R_{\phi}+\frac{\tau-1}{\tau}|\nabla f|^{2}-(m-\tau)\lambda=\varrho e^{\frac{2}{\tau}f}$$
for some real constant $\varrho$. By the compactness of manifold $M$, there exists some global maximum point $p\in M$ of the function $\varrho e^{\frac{2}{\tau}f}$. Then, for any
point $x\in M$, we have
$$\aligned
R_{\phi}(p)-(m-\tau)\lambda&=\left(\varrho e^{\frac{2}{\tau}f}\right)(p)\\
&\geq\left(\varrho e^{\frac{2}{\tau}f}\right)(x)=R_{\phi}(x)+\frac{\tau-1}{\tau}|\nabla f|^{2}(x)-(m-\tau)\lambda\\
&\geq R_{\phi}(x)-(m-\tau)\lambda,
\endaligned$$
which implies that $R_{\phi}$ achieves its maximal value $(R_{\phi})_{max}$ at $p$ and
$$\frac{\tau-1}{\tau}|\nabla f|^{2}\leq (R_{\phi})_{max}-R_{\phi}$$
holds for all $x\in M$, and we obtain \eqref{lem5-1-3}.
\endproof

\subsection{Proof of Theorems \ref{thm3}}
It was proved in \cite{Wanglf} that when $\tau> 1$,
$$R_{\phi}\geq\frac{m(m-1)}{m+\tau-1}\lambda.$$
Therefore, by \eqref{lem5-1-2} we have
\begin{equation}\label{thm3-pf1}\aligned
&\int\limits_{M^m}(\Delta_{\tau,f}(f))^{2}\,e^{-\frac{f}{\tau}}d\upsilon\\
\geq&\int\limits_{M^m}\left(\left(2+m-\frac{(3\tau+m-1)}{\tau(m+\tau-1)}\right)\lambda|\nabla f|^{2}-\frac{\tau-1}{\tau}R_{\phi}|\nabla f|^{2}\right)\,e^{-\frac{f}{\tau}}d\upsilon\\
\geq&\int\limits_{M^m}\left((\tau-1)\left(\frac{2}{m+\tau-1}+\frac{m}{\tau}\right)\lambda|\nabla f|^{2}+R_{\phi}(R_{\phi}-(R_{\phi})_{max})\right)\,e^{-\frac{f}{\tau}}d\upsilon,
\endaligned\end{equation}
where in the third inequality we have used \eqref{lem5-1-3}.
Let $$V_{\tau,f}=\int\limits_{M^m}\,e^{-\frac{f}{\tau}}d\upsilon.$$
Integrating \eqref{lem2-4-0} on $M$ with respect to the measure $e^{-\frac{f}{\tau}}d\upsilon$ leads to
\begin{equation}\label{thm3-pf3}\int\limits_{M^m}R_{\phi}\,e^{-\frac{f}{\tau}}d\upsilon=m\lambda V_{\tau,f}.\end{equation}
Note that $R_{\phi}=m\lambda-\Delta_{\tau,f}(f)$, which gives
$$R_{\phi}^{2}=(m\lambda)^{2}+(\Delta_{\tau,f}(f))^{2}-2m\lambda\Delta_{\tau,f}(f).$$
Integrating the above equation on $M$ with respect to the measure $e^{-\frac{f}{\tau}}d\upsilon$ leads to
\begin{equation}\label{thm3-pf2}
\int\limits_{M^m}R_{\phi}^{2}\,e^{-\frac{f}{\tau}}d\upsilon=(m\lambda)^{2}V_{\tau,f}+\int\limits_{M^m}(\Delta_{\tau,f}(f))^{2}\,e^{-\frac{f}{\tau}}d\upsilon.
\end{equation}
By \eqref{thm3-pf1}, \eqref{thm3-pf3} and \eqref{thm3-pf2}, we have
$$\aligned
&\int\limits_{M^m}(\Delta_{\tau,f}(f))^{2}\,e^{-\frac{f}{\tau}}d\upsilon\\
\geq&\int\limits_{M^m}(\tau-1)\left(\frac{2}{m+\tau-1}+\frac{m}{\tau}\right)\lambda|\nabla f|^{2}\,e^{-\frac{f}{\tau}}d\upsilon+\int\limits_{M^m}(\Delta_{\tau,f}(f))^{2}\,e^{-\frac{f}{\tau}}d\upsilon\\
&+((m\lambda)^{2}-m\lambda(R_{\phi})_{max})V_{\tau,f},
\endaligned$$
which gives
\begin{equation}\label{thm3-pf5}(\tau-1)\left(\frac{2}{m(m+\tau-1)}+\frac{1}{\tau}\right)\int\limits_{M^m}|\nabla f|^{2}\,e^{-\frac{f}{\tau}}d\upsilon+(m\lambda-(R_{\phi})_{max})V_{\tau,f}\leq0.\end{equation}
Hence, by \eqref{general7} in the theorem, the equality in \eqref{thm3-pf5} must be achieved. This
shows that the equality in \eqref{lem5-1-3} must also attain. Therefore, by \eqref{lem2-4-3} we have
$$(R_{\phi})_{max}-(m-\tau)\lambda=\varrho e^{\frac{2}{\tau}f}.$$
Hence $f$ is a constant and $(M^m,g)$ is harmonic-Einstein.

\subsection{Proof of Theorems \ref{thm4}} Firstly, we establish the following lemma, which plays crucial roles in the proof of Theorems \ref{thm4}.

\begin{lem}\label{lem5-2}
Let $(M^m,g)$ be an $m$-dimensional compact manifold with a metric $g$ satisfying $\tau$-quasi Ricci-harmonic metric equation \eqref{general3} with $\lambda>0$ and $\tau>64m$. If
\begin{equation}\label{lem5-2-1}
Ric_{\phi}\geq(1-\delta)\lambda g
\end{equation}
for some $\delta$ satisfying $0<\delta<\frac{1}{6}$. Then
\begin{equation}\label{lem5-2-2}
R_{\phi}\leq(m+3\tau)\lambda
\end{equation}
and
\begin{equation}\label{lem5-2-3}
\int\limits_{M^m}|\nabla f|^{2}\,e^{-\frac{f}{\tau}}d\upsilon\leq3m\tau\lambda\delta\int\limits_{M^m}\,e^{-\frac{f}{\tau}}d\upsilon.
\end{equation}
\end{lem}

\proof Since
$ Ric_{\phi}\geq(1-\delta)\lambda g$
for some $\delta$ satisfying $0<\delta<\frac{1}{6}$, the generalized scalar
curvature satisfies
\begin{equation}\label{lem5-2-pf1}
R_{\phi}\geq(1-\delta)\lambda m
\end{equation}
and it follows from Myers theorem that
\begin{equation}\label{lem5-2-pf2}
diam(M, g)\leq\pi\sqrt{\frac{m-1}{(1-\delta)\lambda}}.
\end{equation}
Then, by \eqref{lem2-4-3} and \eqref{lem5-2-pf1}, we have
\begin{equation}\label{lem5-2-pf3}\aligned
|\nabla f|^{2}=&\frac{\tau}{\tau-1}\left(\varrho e^{\frac{2}{\tau}f}-(\tau-m)\lambda-R_{\phi}\right)\\
\leq&\frac{\tau\varrho}{\tau-1}e^{\frac{2}{\tau}f}-\frac{(\tau-m\delta)\tau}{\tau-1}\lambda.
\endaligned\end{equation}
\eqref{lem5-2-pf3} can be rewritten as
$$\left|\nabla \arctan\left(\frac{1}{\sqrt{(\tau-m\delta)\lambda}}\sqrt{\varrho e^{\frac{2}{\tau}f}-(\tau-m\delta)\lambda} \right) \right|\leq\sqrt{\frac{(\tau-m\delta)\lambda}{\tau(\tau-1)}}.$$
By the compactness of manifold $M$, there exists some global minimum point $q\in M$ of the potential function $f$. Due to \eqref{lem2-4-4}, the maximum principle shows that
$$\varrho e^{\frac{2}{\tau}f}(q)\leq \tau\lambda,$$
which gives
$$\frac{1}{\sqrt{(\tau-m\delta)\lambda}}\sqrt{\varrho e^{\frac{2}{\tau}f}-(\tau-m\delta)\lambda}\leq\sqrt{\frac{m\delta}{\tau-m\delta}}.$$
Therefore, for any point $x\in M$,
$$\aligned
&\arctan\left(\frac{1}{\sqrt{(\tau-m\delta)\lambda}}\sqrt{\varrho e^{\frac{2}{\tau}f}-(\tau-m\delta)\lambda} \right)(x)\\
&-\arctan\left(\frac{1}{\sqrt{(\tau-m\delta)\lambda}}\sqrt{\varrho e^{\frac{2}{\tau}f}-(\tau-m\delta)\lambda} \right)(q)\\
\leq&\left(\max_{M}\left|\nabla \arctan\left(\frac{1}{\sqrt{(\tau-m\delta)\lambda}}\sqrt{\varrho e^{\frac{2}{\tau}f}-(\tau-m\delta)\lambda} \right) \right|\right)diam(M, g)\\
=&\sqrt{\frac{(\tau-m\delta)\lambda}{\tau(\tau-1)}}diam(M, g)
\leq\pi\sqrt{\frac{(\tau-m\delta)(m-1)}{\tau(\tau-1)(1-\delta)}}.
\endaligned$$
Then, it follows from the above that
\begin{equation}\label{lem5-2-pf5}\aligned
&\arctan\left(\frac{1}{\sqrt{(\tau-m\delta)\lambda}}\sqrt{\varrho e^{\frac{2}{\tau}f}-(\tau-m\delta)\lambda} \right)(x)\\
\leq&\arctan\sqrt{\frac{m\delta}{\tau-m\delta}}+\pi\sqrt{\frac{(\tau-m\delta)(m-1)}{\tau(\tau-1)(1-\delta)}}.
\endaligned\end{equation}
Since $\tau>64m$ we know that
$$\frac{m-1}{\tau-1}\leq\frac{1}{64}$$
and
$$\frac{\tau-m\delta}{\tau(1-\delta)}=\frac{m}{\tau}+\frac{\tau-m}{\tau(1-\delta)}<4.$$
We deduce from \eqref{lem5-2-pf5} that (see \cite{Wanglinf2018})
$$\aligned
\varrho e^{\frac{2}{\tau}f}\leq&(\tau-m\delta)\lambda+\left(\frac{1+\sqrt{\frac{m\delta}{\tau-m\delta}}}{1-\sqrt{\frac{m\delta}{\tau-m\delta}}}\right)^{2}(\tau-m\delta)\lambda\\
=&\frac{2(\tau-m\delta)\lambda\tau}{\tau-2\sqrt{m\delta(\tau-m\delta)}}\leq4\tau\lambda.
\endaligned$$
By \eqref{lem2-4-3}, we have $$R_{\phi}\leq\varrho e^{\frac{2}{\tau}f}-(\tau-m)\lambda\leq(m+3\tau)\lambda.$$
Hence, we arrive at \eqref{lem5-2-2}.

Now, we prove \eqref{lem5-2-3}. For simplicity, we put
$$\Omega_{+}:=\{x\in M: R_{\phi}(x)>m\lambda\}\,\,\,\,\,\,\,\,\,and\,\,\,\,\,\,\,\,\,  \Omega_{-}:=\{x\in M: R_{\phi}(x)<m\lambda\},$$ respectively.
By \eqref{lem5-2-pf1} and \eqref{lem5-2-2} we know that
$$(1-\delta)\lambda m\leq R_{\phi}\leq(m+3\tau)\lambda.$$
Then, we have
$$\aligned
&\frac{1}{V_{\tau,f}}\int\limits_{M^m}(\Delta_{\tau,f}(f))^{2}\,e^{-\frac{f}{\tau}}d\upsilon=\frac{1}{V_{\tau,f}}\int\limits_{M^m}(R_{\phi}-m\lambda)^{2}\,e^{-\frac{f}{\tau}}d\upsilon\\
=&\frac{1}{V_{\tau,f}}\int\limits_{\Omega_{+}}(R_{\phi}-m\lambda)^{2}\,e^{-\frac{f}{\tau}}d\upsilon+\frac{1}{V_{\tau,f}}\int\limits_{\Omega_{-}}(R_{\phi}-m\lambda)^{2}\,e^{-\frac{f}{\tau}}d\upsilon\\
\leq&\frac{3\tau \lambda}{V_{\tau,f}}\int\limits_{\Omega_{+}}(R_{\phi}-m\lambda)\,e^{-\frac{f}{\tau}}d\upsilon+m^{2}\delta^{2}\lambda^{2}.
\endaligned$$
Due to \eqref{lem2-4-0} we have
$$\int\limits_{\Omega_{+}}(R_{\phi}-m\lambda)\,e^{-\frac{f}{\tau}}d\upsilon+\int\limits_{\Omega_{-}}(R_{\phi}-m\lambda)\,e^{-\frac{f}{\tau}}d\upsilon=0.$$
We then get that
\begin{equation}\label{lem5-2-pf7}\aligned
\frac{1}{V_{\tau,f}}\int\limits_{M^m}(\Delta_{\tau,f}(f))^{2}\,e^{-\frac{f}{\tau}}d\upsilon\leq&\frac{3\tau \lambda}{V_{\tau,f}}\int\limits_{\Omega_{-}}(m\lambda-R_{\phi})\,e^{-\frac{f}{\tau}}d\upsilon+m^{2}\delta^{2}\lambda^{2}\\
\leq&m\delta(3\tau+m\delta)\lambda^{2}.
\endaligned\end{equation}
Due to \eqref{lem5-1-1} and \eqref{lem5-2-pf1} we have
\begin{equation}\label{lem5-2-pf8}\aligned
&\int\limits_{M^m}(\Delta_{\tau,f}(f))^{2}\,e^{-\frac{f}{\tau}}d\upsilon\\
=&\int\limits_{M^m}\left(\frac{2(\tau-1)}{\tau}Ric_{\phi}(\nabla f,\nabla f)+\frac{2}{\tau}[R_{\phi}-(m-1)\lambda]|\nabla f|^{2}\right)\,e^{-\frac{f}{\tau}}d\upsilon\\
\geq&\int\limits_{M^m}\left(\frac{2(\tau-1)}{\tau}Ric_{\phi}(\nabla f,\nabla f)+\frac{2\lambda}{\tau}(1-m\delta)|\nabla f|^{2}\right)\,e^{-\frac{f}{\tau}}d\upsilon.
\endaligned\end{equation}
Therefore, by \eqref{lem5-2-pf7} and \eqref{lem5-2-pf8} we have
$$\aligned
&\frac{2(\tau-1)}{\tau V_{\tau,f}}\int\limits_{M^m}Ric_{\phi}(\nabla f,\nabla f)\,e^{-\frac{f}{\tau}}d\upsilon+\frac{2\lambda}{\tau V_{\tau,f}}\int\limits_{M^m}|\nabla f|^{2}\,e^{-\frac{f}{\tau}}d\upsilon\\
\leq&m\delta(3\tau+m\delta)\lambda^{2}+\frac{2m\delta\lambda}{\tau V_{\tau,f}}\int\limits_{M^m}|\nabla f|^{2}\,e^{-\frac{f}{\tau}}d\upsilon.
\endaligned$$
By \eqref{lem5-2-1}, we get that
$$(2\tau+2(1-m-\tau)\delta)\int\limits_{M^m}|\nabla f|^{2}\,e^{-\frac{f}{\tau}}d\upsilon\leq m\tau\delta\lambda(3\tau+m\delta)V_{\tau,f}.$$
By the fact that $\tau>64m$ and $0<\delta<\frac{1}{6}$, it is easy to verify that
$$(2\tau+2(1-m-\tau)\delta)\geq\frac{1}{3}(3\tau+m\delta).$$
Hence, we arrive at \eqref{lem5-2-3}.
\endproof
Now, we may finish the proof of Theorem \ref{thm4}. We here assume that $(M^m,g)$ is not harmonic-Einstein and deduce a contradiction.
If $(M^m,g)$ is not harmonic-Einstein, by Lemma \ref{lem5-2}, we know that
$$\int\limits_{M^m}|\nabla f|^{2}\,e^{-\frac{f}{\tau}}d\upsilon\leq3m\tau\lambda\delta\int\limits_{M^m}\,e^{-\frac{f}{\tau}}d\upsilon,$$
which contradicts \eqref{general9}. This completes the proof of Theorem \ref{thm4}.

\bibliographystyle{Plain}

\end{document}